\documentclass[12pt]{article}

\usepackage{amsmath,amsfonts}
\usepackage{multicol}
\usepackage{graphicx}

\newcommand{\CC}{{\cal{C}}}

\newcommand{\LL}{{\cal{L}}}
\newcommand{\R}{\mathbb{R}}
\newcommand{\N}{\mathbb{N}}
\renewcommand{\SS}{{\cal{S}}}

\def\z{\mathbf{z}}
\def\w{\mathbf{w}}
\def\y{\mathbf{y}}
\def\K{\mathbf{K}}

\newtheorem{assumption}{Assumption}
\newtheorem{example}{Example}
\newtheorem{remark}{Remark}

\title{Nonlinear optimal control synthesis\\ via occupation measures}
\author{Didier Henrion, Jean B. Lasserre and Carlo Savorgnan
\thanks{The authors are with LAAS-CNRS, University of Toulouse, France.
D. Henrion is also with the Faculty of Electrical Engineering, Czech Technical University in Prague, Czech Republic. J. B. Lasserre is also with the Institute of Mathematics of the University
of Toulouse, France. E-mail: {\tt \{henrion,lasserre,savorgnan\}@laas.fr}}}

\begin{document}
\maketitle

\begin{abstract}
We consider nonlinear optimal control problems (OCPs) for which all problem data are polynomial. In the first part of the paper, we review how occupation measures can be used to approximate pointwise the optimal value function of a given OCP, using a hierarchy of linear matrix inequality (LMI) relaxations. In the second part, we extend the methodology to approximate the optimal value function on a given set and we use such a function to constructively and computationally derive an almost optimal control law.
Numerical examples show the effectiveness of the approach.
\end{abstract}

\section{Introduction}

It is well known that solving an optimal control problem (OCP) can be a very hard task notwithstanding the power of theoretical tools such Pontryagin's minimum principle and Hamilton-Jacobi-Bellman optimality condition.
This statement is particularly true when dealing with state and input constraints.

\textbf{Contribution.} In this paper we consider the class of OCPs for which all problem data are polynomial.
The approach we deploy (which was introduced in \cite{LasPriHen2005}) is based on moment theory and consists in deriving a hierarchy of convex linear matrix inequality (LMI) relaxations of the OCP which give an increasing sequence of lower bounds on the optimal value. These LMI problems can be solved using off-the-shelf semidefinite programming (SDP) solvers.

The contribution with respect to \cite{LasPriHen2005} and its extended version \cite{LasHenPriTre2008} is twofold. First, the derivation of the relaxation is obtained in a simpler way, starting from basic concepts. The second and more important contribution is that we show how the methodology can be applied to approximate the optimal value function on a set and to derive constructively and computationally
a control law. The approach is illustrated on a few simple examples.

\textbf{Notation.}
$\R$ and $\N$ denote respectively the sets of real and integer numbers.
$\R[y]=[y_1, \dots , y_n]$ denotes the ring of polynomials in the variable $y$.
$\R[y]_d=[y_1, \dots , y_n]$ denotes the ring of polynomials of degree at most $d$ in the variable $y$.
When $y \in \R^n$ and $\alpha\in\N^n$, $y^\alpha$ stands for $y_1^{\alpha_1} \dots y_n^{\alpha_n}$.
Given a polynomial function $\varphi$, $\deg(\varphi)$ is the maximal degree of its monomials.
Given a differentiable function $\varphi(y)$, $\nabla_y(\varphi)=[\frac{\partial\varphi}{\partial y_1}, \dots ,\frac{\partial\varphi}{\partial y_n}]$ is its gradient with respect to $y$.
$\delta_{y_0}$ is the Dirac measure at $y_0$.
$v'$ denotes the transpose of $v$.

\section{Problem definition}
Consider continuous-time systems described by the differential equation
\begin{equation}\label{eq:dynamics}
\dot x(t) = f(t, x(t), u(t))
\end{equation}
where $x\in\R^n$ and $u\in\R^m$ are respectively the state vector and input vector.
By defining the cost function
\begin{equation}\label{eq:cost}
\int_0^T h(t, x(t), u(t)) dt + H(x(T))
\end{equation}
the initial constraint
\begin{equation}\nonumber
x(0) \in \CC_I=\{ x: g_{I_j}(x)\leq 0, ~ j=1,\dots,n_I \} 
\end{equation}
and the final constraint
\begin{equation}\nonumber
x(T) \in \CC_F=\{ x: g_{F_j}(x)\leq 0, ~ j=1,\dots,n_F \}
\end{equation}
we can formulate several OCPs. E.g., when $\CC_I$ and $\CC_F$ contain only one point we have the classical problem of driving the system from an assigned initial condition $x(0)=x_0$ to a final condition $x(T)=x_T$ by minimizing a given cost.

In the sequel, we will consider all the problems which can be cast in this framework with the additional constraint on the trajectory $(t, x(t), u(t)) \in \CC_T$, where
\begin{equation}\nonumber
\CC_T=\{ (t,x,u): g_{T_j}(t,x,u)\leq 0, ~ j=1,\dots,n_T \}.
\end{equation}

An important assumption which is necessary for the derivation of the methodology is that all problem data are polynomial. More precisely:
\begin{assumption}\label{as:poly}
The functions $f$, $h$, $H$, $g_{I_j}$, $g_{T_j}$ and $g_{F_j}$ are polynomial.
\end{assumption}

\section{The moment approach to optimal control}
The key idea underlying the moment approach is that of defining three \textit{occupation measures} which convey the information about the initial condition of the system, its trajectory and the final condition.
The OCP is then rephrased in terms of the moments of such measures.
The convex problem obtained contains three ingredients:
\begin{itemize}
\item a set of linear equality constraints on the moments which characterize the system dynamics;
\item a set of semidefinite constraints which come from the fact that the moments belong to a measure;
\item a set of semidefinite constraints which translate the constraints induced by $\CC_I$, $\CC_T$ and $\CC_F$ on the supports of the measures.
\end{itemize}
To derive this constraints we assume horizon $T$ is fixed.

\subsection{The trajectory constraints}
To obtain the trajectory constraints we start from the idea that the system trajectories can be characterized studying how certain \textit{test functions} evolve along the trajectories. For this purpose, we choose functions which are monomials of the form $t^\alpha x^\beta$. Consider a trajectory $x(t)$. Using the fundamental theorem of calculus we can write
\begin{equation}\label{eq:ftoc}
T^\alpha x(T)^\beta = 0^\alpha x(0)^\beta + \int_0^T \frac{d(t^\alpha x(t)^\beta)}{dt} dt.
\end{equation}
The trajectory constraints are obtained by rephrasing equation (\ref{eq:ftoc}) in terms of 
three properly defined occupation measures.

The \textit{final occupation measure} $\mu_F$ captures the information on the state at time $T$
\begin{equation}\nonumber
x(T)^\beta = \int x^\beta \delta_{x(T)}(dx) = \int x^\beta d\mu_F.
\end{equation}
The \textit{initial occupation measure} $\mu_I$ captures the information on the initial condition of the system
\begin{equation}\nonumber
x(0)^\beta = \int x^\beta \delta_{x(0)}(dx) = \int x^\beta d\mu_I.
\end{equation}
The \textit{trajectory occupation measure} $\mu_T$ captures the information on the value of $t$, $x(t)$ and $u(t)$ along the trajectory
\[
\int_0^T \!\! t^\gamma x(t)^\eta u(t)^\nu dt \!
=\!\int_0^T \!\!\!\!\int t^\gamma x^\eta u^\nu \delta_{x(t),u(t)}(dx,du)dt \!
=\!\int x^\gamma u^\eta t^\nu d\mu_T.
\]
Notice that $\mu_F$ and $\mu_I$ are probability measures, their mass is equal to $1$.

Next, if $f\in\R[t,x,u]$ then $\frac{d(x^\alpha t^\beta)}{dt}\in\R[t,x,u]$ because
\begin{equation}\nonumber 
\frac{d(x^\alpha t^\beta)}{dt} \! = \! \frac{\partial(x^\alpha t^\beta)}{\partial t} + \nabla_x (x^\alpha t^\beta)f(x,u) \! = \!\! \sum_{\gamma,\eta,\nu} \! a^{\alpha\beta}_{\gamma\eta\nu} x^\gamma u^\eta t^\nu
\end{equation}
for some coefficients $a^{\alpha\beta}_{\gamma\eta\nu}$ that depend on $f$.
The degree of the derivative is $\deg(x^\alpha t^\beta)-1+\deg(f)$.
Using previous equations, (\ref{eq:ftoc}) can be rephrased as
\begin{equation}\label{eq:meascons}
T^\alpha \int x^\beta d\mu_F 
= 0^\alpha \int x^\beta d\mu_I + \sum_{\gamma,\eta,\nu} a^{\alpha\beta}_{\gamma\eta\nu} \int t^\gamma x^\eta u^\nu d\mu_T,
\end{equation}
i.e., a {\it linear} relationship between the moments of 
$\mu_F,\mu_I$ and $\mu_T$.
Namely, introducing the notation 
$z_{\beta}=\int x^\beta d\mu_F$,
$w_{\beta}=\int x^\beta d\mu_I$,
$y_{\gamma\eta\nu}=\int t^\gamma x^\eta u^\nu d\mu_T$, 
we obtain
\begin{equation}\label{eq:momcons}
T^\alpha z_{\beta} = 0^\alpha w_{\beta} + \sum_{\gamma,\eta,\nu} a^{\alpha\beta}_{\gamma\eta\nu} y_{\gamma\eta\nu}
\end{equation}
for every $\alpha,\beta\in\N\times\N^{n}$. Notice that from (\ref{eq:momcons}),
the mass of $\mu_T$ is $T$.
In compact notation, consider test functions of degree up to $r$ and the canonical basis of monomials of degree at most $r$:
\begin{equation}\nonumber
m_r(x)=[1,x_1,\dots,x_n,x_1^2,x_1x_2,\dots,x_1^r,x_1^{r-1}x_2,\dots,x_n^r]'.
\end{equation}
Define the vectors 
$\z_r=\int m_r(x) d\mu_F$,
$\w_r=\int m_r(x) d\mu_I$ and
$\y_k=\int m_k(t,x,u) d\mu_T$.
Then,
\begin{equation}\label{eq:momcon}
A_F \z_r = A_I \w_r + A_T \y_k
\end{equation}
where $k\geq r-1+\deg(f)$ and the coefficients of the matrices $A_F$, $A_I$ and $A_T$ can be obtained from equation (\ref{eq:momcons}).

Define the coefficient vectors $c_h$ and $c_H$ to be such that
\begin{equation}\nonumber
h(t,x,u)=c_h'm_k(x,t,u), \qquad H(x)=c_H' m_r(x).
\end{equation}
Observe that
\begin{equation}\label{criterion}
\int_0^T h(t, x(t), u(t)) dt + H(x(T))=c_h' \y_k+  c_H' \z_r,
\end{equation}
i.e., the criterion of the OCP is a linear functional on $\z_r$ and $\y_k$.

So far, for a given trajectory $x(t)$, we have characterized
linear constraints satisfied by the moments of the three
associated occupation meaures. Now, if the trajectory is unknown,
the three measures are unknown, and we can consider the abstract
linear programming (LP) problem $J(\mu_I) = \min_{\mu_T,\mu_I,\mu_F}
\int h d\mu_T + \int H d\mu_F$
subject to (\ref{eq:meascons}),
which aims at finding the occupation measures associated with
the optimal trajectory. The measures $\mu_F$, $\mu_I$, $\mu_T$
are characterized through their respective
truncated moment vectors $\z_r$, $\w_r$, $\y_k$,
the remaining difficulty being finding conditions that ensure
that those vectors and indeed moment vectors of measures
with respective supports $\CC_F$, $\CC_I$, $\CC_T$. This is explained in the
next section.

A nice feature of the approach is that we can play with the initial
and final measures. For instance, if $\mu_I=\delta_{x_0}$
we retrieve the optimal cost $J(\delta_{x_0})$ of the OCP with fixed
initial state $x_0$. Now, if $\mu_I$ is unknown, but with
known support $\CC_I$, then $J(\mu_I)=\min_{x_0\in\CC_I}J(\delta_{x_0})$.
Finally, if $\mu_I$ is known, but not a Dirac, solving the
above LP problem aims at computing $\int J(\delta_{x_0})d\mu_I(x_0)$.

\subsection{The moment matrix constraints}
There exist linear programming (LP) or semidefinite programming (SDP)
necessary and sufficient conditions 
for an infinite vector to be a {\it moment} vector, i.e., the vector of moments
of some finite Borel measure on a compact 
basic semi-algebraic set; see e.g. \cite{Put1993}.
We chose the latter since it has shown to be more effective for numerical purposes \cite{LasPri2004}.

With $r$ an even number, let
\begin{equation}\nonumber
M(\z_r)=\int m_{r/2}(x) m_{r/2}(x)' d\mu_F
\end{equation}
be the {\it moment matrix} of order $r$ associated with $\mu_F$.
Obviously, $M(\z_r)$ is positive semidefinite, denoted $M(\z_r)\succeq0$. Therefore,
in the convex relaxation of the OCP, one imposes
\begin{equation}\label{eq:mmc}
M(\z_r) \succeq 0,
\end{equation}
and similar constraints are imposed on $\w_r$ and $\y_k$.

\subsection{The localizing matrix constraints}
Similarly to the previous subsection, one may express the support constraints induced by $\CC_I$, $\CC_T$ and $\CC_F$ in terms of linear matrix inequalities on $\z_r,\w_r$ and $\y_k$. To derive such inequalities, define
\begin{equation}\nonumber
d_{F_j}=\left\{
\begin{array}{ll}
\deg(g_{F_j}(x)) & \mbox{if}~\deg(g_{F_j}(x))~\mbox{is even} \\
\deg(g_{F_j}(x))+1 & \mbox{if}~\deg(g_{F_j}(x))~\mbox{is odd}
\end{array} \right.
\end{equation}
and the localizing matrix
\begin{equation}\nonumber
L_{g_{F_j}}(\z_r)=\int g_{F_j}(x) m_{(r-d_{F_j})/2}(x) m_{(r-d_{F_j})/2}(x)' d\mu_F.
\end{equation}
 The matrix $g_{F_j}(x) m_{(r-d_{F_j})/2}(x) m_{(r-d_{F_j})/2}(x)'$ is positive semidefinite for every value of $x$ such that $g_{F_j}(x)\geq 0$. Hence if $\mu_F$ is supported  on $\CC_F$ then
 $L_{{F_j}}(\z_r)\succeq0$ for every $j$. Therefore, in the convex relaxation of the OCP one imposes
the semidefinite constraint
\begin{equation}\label{eq:lmc}
L_{g_{F_j}}(\z_r) \succeq 0 \qquad j=1,\dots,n_F
\end{equation}
and similar  semidefinite constraints for $\w_r$ and $\y_k$.

Further details on moment and localizing matrix constraints can be found in \cite{Las2001}.

\subsection{The convex relaxation}
To construct the convex relaxation of the OCP, let $r$ and $k$ be even numbers such that
\begin{equation}\nonumber
r\geq\deg(H), \quad k\geq\deg(h), \quad k\geq r+\deg(f).
\end{equation}
In this paper, we will assume that the initial probability
measure $\mu_I$ is known through its moments $\w_r$. 

The convex relaxation is the following truncated moment problem:
\begin{equation}\label{eq:momocp2}
\begin{split}
\min_{\z_r, \y_k}\quad & c_h' \y_k+  c_H' \z_r \\
& A_F \z_r = A_I \bar \w_r + A_T \y_k \\
& M(\z_r)\succeq 0, ~~ L_{g_{F_j}}(\z_r)\succeq 0, ~~ \forall j=1,\dots,n_F \\
& M(\y_k)\succeq 0, ~~ L_{g_{T_j}}(\y_k)\succeq 0, ~~ \forall j=1,\dots,n_T
\end{split}
\end{equation}
where the notation $\bar \w_r$ indicates that the moment vector is known.
Two important facts should be noticed for the moment problem (\ref{eq:momocp2}):
\begin{itemize}
\item the constraints on the moments correspond to necessary conditions and therefore, 
in general one only obtains a lower bound on the optimal value of OCP;
\item with $\hat r > r$ and $\hat k > k$, the constraints of the original problem with $r$ and $k$ are a subset of the constraints of the problem with $\hat r$ and $\hat k$. Therefore, increasing the value of $r$ and $k$ yields a monotonically nondecreasing sequence of lower bounds on the optimal value.
\end{itemize}

\begin{remark}
If the initial measure $\mu_I$ was unknown, we would have to include the additional
constraints $M(\w_r)\succeq 0$, $L_{g_{I_j}}(\w_r)\succeq 0$, $\forall j=1,\dots,n_I$
with now $\w_r$ being an unknown moment vector with first entry equal to one.
\end{remark}

\begin{remark}
One goal of this paper is to derive the convex relaxation of the OCP starting from really basic notions. The same optimization problem can be also obtained using as a starting point the duality between the Banach space of bounded continuous functions on a compact set $\K$ and the Banach space of finite signed Borel measures on $\K$, as done in \cite{LasPriHen2005,LasHenPriTre2008} where the sequence of lower bounds on the optimal value is shown to converge under some assumptions on the problem data. The interested reader is referred to these paper for further details.
\end{remark}

\section{The dual approach: SOS polynomials}

For the developments of the results in sections \ref{se:set_approx} and \ref{se:control} and a better understanding of the moment approach to OCP, it is important to look at its dual formulation which has an interesting interpretation in terms of SOS polynomials.

A polynomial $p\in\R[x]$ of degree $2d$ is an SOS if 
$p(x)=\sum_{i=1}^s f_i(x)^2$ for some $(f_i)_{i=1}^s\in\R[x]$, and this
implies that $p$ is non-negative.
And $p$ is an SOS if and only if there exists a positive semidefinite matrix $Q$ such that $p(x)\equiv m_{d}(x)'Q m_{d}(x)$. Denote by $\Sigma[x]$ the set of SOS polynomials and by $\Sigma[x]_r$ (with $r$ even) the set of SOS polynomials of degree at most $r$. See \cite{Par2003} for more details.

The SDP dual of  (\ref{eq:momocp2}) is
\begin{equation}\label{eq:dualocpsdp}
\begin{split}
\max_{\substack{c_\varphi, S\succeq 0, Q\succeq 0 \\ S_j\succeq 0, Q_j\succeq 0}} ~ & (A_I \bar \w_r )' c_\varphi \\
& -A_T' c_\varphi + M^*(S) + \sum_{j=1}^{n_T} L_{g_{T_j}}^*(S_j) = c_h \\
& A_F' c_\varphi + M^*(Q) + \sum_{j=1}^{n_T} L_{g_{F_j}}^*(S_j) =c_H
\end{split}
\end{equation}
where the symbol $*$ indicates the adjoint operator, $c_\varphi$ is the dual variable associated with
the (moment) trajectory constraint, $S$ and $Q$ are the dual variables associated with the moment matrix constraints and $S_j$ and $Q_j$ are the dual variables associated with the localizing matrix constraints.
To interpret problem (\ref{eq:dualocpsdp}) in terms of SOS polynomials, define the polynomial function
\begin{equation}\nonumber
(t,x)\longmapsto \varphi(t,x)=c_\varphi m_r(t,x).
\end{equation}
Explicitating the adjoint operators yields the following problem:
\begin{equation}\label{eq:sosocp}
\begin{split}
\max_{\substack{\varphi\in\R[t,x]_r, s\in \Sigma[t,x,u]_k, q\in \Sigma[x]_r \\ s_j\in \Sigma[t,x,u]_{k-d_{T_j}}, q_j\in
\Sigma[x]_{r-d_{F_j}}}} ~ & \varphi(0,x(0)) \\
& \hspace{-3cm}
\frac{\partial\varphi(x,t)}{\partial t} + \nabla_x \varphi(x,t)f(t,x,u) + h(t,x,u) = s(t,x,u)+ \sum_{j=1}^{n_T} g_{T_j}(t,x,u)s_j(t,x,u) \\
& \hspace{-3cm}
\varphi(x,T) - H(x) = -q(x)-\sum_{j=1}^{n_F} g_{F_j}(x)q_j(x).
\end{split}
\end{equation}
Consider the right handsides of the first and second constraints in (\ref{eq:sosocp}). The first one is a polynomial non-negative on $\CC_T$, while the second one is a polynomial non-positive on $\CC_F$.
In fact, both are Putinar's SOS representations of 
their respective left-hand-side \cite{Put1993}.
 As a consequence, every feasible solution $\varphi$ of (\ref{eq:sosocp}) is such that
\begin{equation}\label{eq:postraj}
\frac{\partial\varphi(x,t)}{\partial t} + \nabla_x \varphi(x,t) f(t,x,u) + h(t,x,u) \geq 0
\end{equation}
for all $(t,x,u) \in \CC_T$ and
\begin{equation}\label{eq:posfin}
H(x) - \varphi(T,x) \geq 0 \qquad \forall x \in \CC_F.
\end{equation}
Suppose the OCP has a solution and consider an optimal control law $\bar u(t)$ which generates an optimal trajectory $\bar x(t)$. Therefore, the optimal value function is
\[
\bar \varphi(t,\bar x(t))=\int_0^T h(t, \bar x(t), \bar u(t)) dt + H(\bar x(T)).
\]
Since a solution $\varphi(t,x)\in\R[t,x]$ of the optimization problem (\ref{eq:sosocp})
is differentiable, the fundamental theorem of calculus yields
\begin{multline}\nonumber
\varphi(t,\bar x(t))=\varphi(T, \bar x(T))-\int_t^T \frac{\varphi(\bar x(\theta),\theta)}{\partial \theta} + \\
\nabla_x \varphi(\bar x(\theta),\theta)f(\theta, \bar x(\theta),  \bar u(\theta)) d\theta .
\end{multline}
Combining with the two preceding equations yields
\begin{equation}\label{eq:intdiff}
\begin{array}{l}
\bar \varphi(t,\bar x(t)) \! - \! \varphi(t,\bar x(t)) \! = \!
[H(\bar x(T)) \! - \! \varphi(T,\bar x(T))] \! + \\
\quad [\int_t^T\frac{\partial\varphi(\bar x(\theta),\theta)}{\partial \theta} \! + \! \nabla_x \varphi(\bar x(\theta),
\theta) f(\theta, \bar x(\theta),  \bar u(\theta)) \! + \! h(\theta, \bar x(\theta), \bar u(\theta)) d\theta].
\end{array}
\end{equation}
From (\ref{eq:postraj}) and (\ref{eq:posfin}), both terms in square braces
in the right handside of (\ref{eq:intdiff}) are non-negative. Therefore:

\begin{itemize}
\item $\bar \varphi(0,x(0)) - \varphi(0, x(0))\geq0$ and so, as in the moment formulation, one obtains a lower bound on the optimal value function.
Therefore, the SOS formulation can be interpreted as the search of a smooth subsolution of the Hamilton-Jacobi-Bellman optimality condition
\begin{equation}\nonumber
\min_{u\in U(t,x)}\left[\frac{\partial\varphi(x,t)}{\partial t} + \nabla_x \varphi(x,t) + h(t,x,u)\right] = 0
\end{equation}
($U(t,x)$ being the set of admissible input at $x(t)=x$).
\item if $\bar \varphi(0,x(0)) - \varphi(0, x(0))$ is small, both terms in the right handside of
(\ref{eq:intdiff}) are small. This implies that the integrand of the second term
is small all along the trajectory.
\end{itemize}

\begin{remark}\label{re:no_t}
So far we have discussed only OCPs for which the value of $T$ is fixed. The moment approach 
also applies for problems with free terminal time and when the dynamics
does not depend on $t$. By Bellman's principle of optimality, the optimal value function does not depend on $t$ and so, in this case, the test functions are of the form $x^\beta$ and the occupation measure $\mu_T$ is on $\R^n\times\R^m$ (instead of $\R\times\R^n\times\R^m)$.
\end{remark}

We next illustrate the effectiveness of the approach 
on a simple numerical example which also motivates the developments of the next section.
\begin{example}\label{ex:dint}
Consider the double integrator 
\begin{equation}\nonumber
\begin{bmatrix} \dot x_1(t) \\ \dot x_2(t) \end{bmatrix} =
\begin{bmatrix} 0 & 1 \\ 0 & 0 \end{bmatrix}
\begin{bmatrix} x_1(t) \\ x_2(t) \end{bmatrix} +
\begin{bmatrix} 0 \\ 1 \end{bmatrix} u(t)
\end{equation}
with constraints $x_2(t)\geq -1$ and  $-1\leq u(t) \leq 1$, $\forall t$.
Driving in minimum time $T(x)$ any initial condition $x(0)=x$ to the origin is an interesting test problem because an analytic solution is available \cite{LasHenPriTre2008}.
Indeed,  with $x_2\geq -1$,
\begin{equation}\nonumber
T(x)=\left\{
\begin{array}{ll}
x_2^2/2+x_1+x_2+1 & \mbox{if}~ x_1 \geq 1-x_2^2/2 \\
2\sqrt{x_2^2/2+x_1}+x_2 & \mbox{if}~  x_1 \leq 1-x_2^2/2 ~\mbox{and} \\
 & x_1 \geq -x_2^2/2 \operatorname{sign}(x_2) \\
2\sqrt{x_2^2/2-x_1}-x_2 & \mbox{if}~ x_1 < -x_2^2/2 \operatorname{sign}(x_2).
\end{array} \right.
\end{equation}
Consider the initial condition $x=(-0.5,-0.8)$ for which
$T(x)=2.6111$. To apply the moment approach, set $h(x(t),u(t))=1$ and $H(x(T))=0$. Solving the moment problem for different values of the degree $r$ yields the values
\vskip 0.2cm
\begin{center}
\begin{tabular}{|l|l|l|l|l|}
\hline degree & 6 & 10 & 14 & 18\\
\hline cost & 1.3882 & 2.1533 & 2.5335 & 2.6061 \\ \hline
\end{tabular}
\end{center}
\vskip 0.2cm
In figure \ref{fig:dintegrator} the function $T(x)$ and $\varphi(x)$ obtained by solving the dual problem for $r=18$ are represented with solid and dashed lines, respectively.
\begin{figure}[htb]
\centering{
\includegraphics[width=0.8\textwidth]{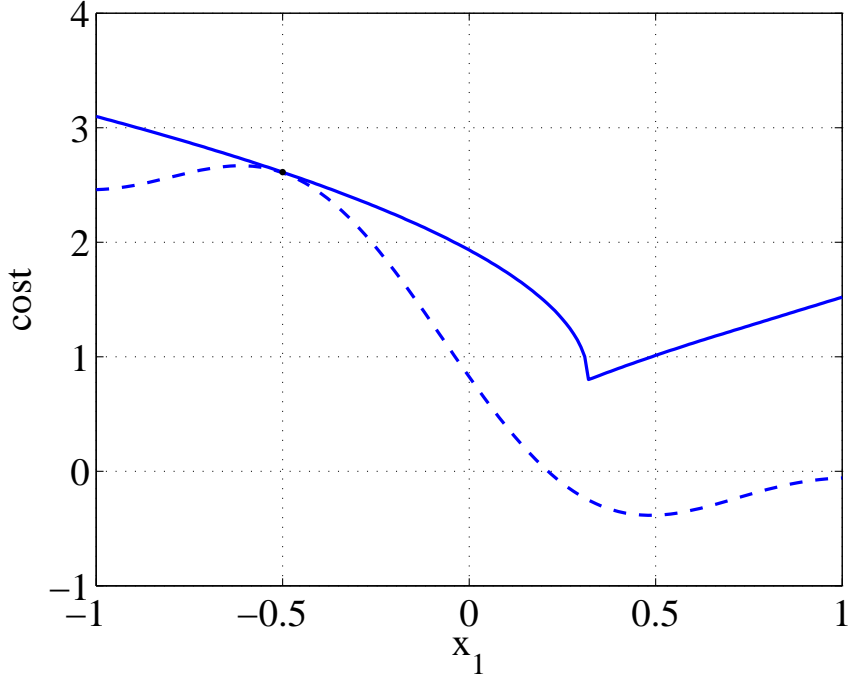} }
\caption{$T(x)$ (solid line) and $\varphi(x)$ (dashed line) from example \ref{ex:dint}. The graph is obtained for $x_2=-0.8$ and $r=18$.}
\label{fig:dintegrator}
\end{figure}
Notice that $\varphi(x)$ approximates very well the value of $T(x)$ at $x=(-0.5,-0.8)$ but it gives a loose lower bound on the other points. As expected, it also gives a very good lower bound on
$T(x(t))$ at every point $x(t)$ of an optimal trajectory from $x(0)=x$.
\end{example}
In the next section, the moment approach to optimal control is extended to obtain good approximations of the value function on a larger set which will be used later to obtain a good control law from the knowledge of an optimal solution $\varphi$ of (\ref{eq:sosocp}).

\section{Approximation of the value function on a set}\label{se:set_approx}
Consider the dual problem (\ref{eq:sosocp}). From the cost $\varphi(0,x(0))$ to be maximized,
the optimal value function $\varphi(t,x)$ is guaranteed to be a good approximation of $T(x)$ 
only at all points $(t,x(t))$ of an optimal trajectory from $x(0)=x$. 
Ideally, for control synthesis purposes, we would like to enlarge the set where $\varphi$ is a good approximation to a region that contains an optimal trajectory $(t,x(t))$ from $x(0)=x$.
To do this, the key observation  is that the trajectory moment constraint (\ref{eq:momcon}) is valid for any initial occupation measure $\mu_I$  and not only for $\mu_I:=\delta_{x}$.

Indeed, let $\mu_I$ be a probability measure for which we know how to calculate its moment vector . The following equality holds
\begin{equation}\nonumber
(A_I \bar w_r )' c_\varphi = \int \varphi(0,x) d\mu_I.
\end{equation}
For instance, if $\mu_I$ is a uniform probability measure on $\SS$
\begin{equation}\nonumber
\int \varphi(0,x) d\mu_I = \int_\SS \varphi(0,x) dx.
\end{equation}
In this case, the solution $\varphi(t,x)$ of the optimization problem minimizes the $\LL_1$ norm $\int_\SS |\bar\varphi(0,x)-\varphi(0,x)| dx$ where $\bar\varphi(t,x)$ is the optimal value function for the problem. This can be easily verified by using the fact that $\int_\SS |\bar\varphi(0,x)-\varphi(0,x)| dx=\int_\SS \bar\varphi(0,x)dx-\int_\SS\varphi(0,x) dx$.
\begin{example}\label{ex:dint_continued}
Consider Example \ref{ex:dint}. Figure \ref{fig:dint_continued} displays  $T(x)$ (solid line) and $\varphi(x)$ when $\SS$ is the line segment $[(-1,-0.8)-(-0.5,-0.8)]$ (dotted line) and when $\SS$ is the segment $[(-1,-0.8),(1,-0.8)]$ (dashed line).
\begin{figure}[htb]
\centering{
\includegraphics[width=0.8\textwidth]{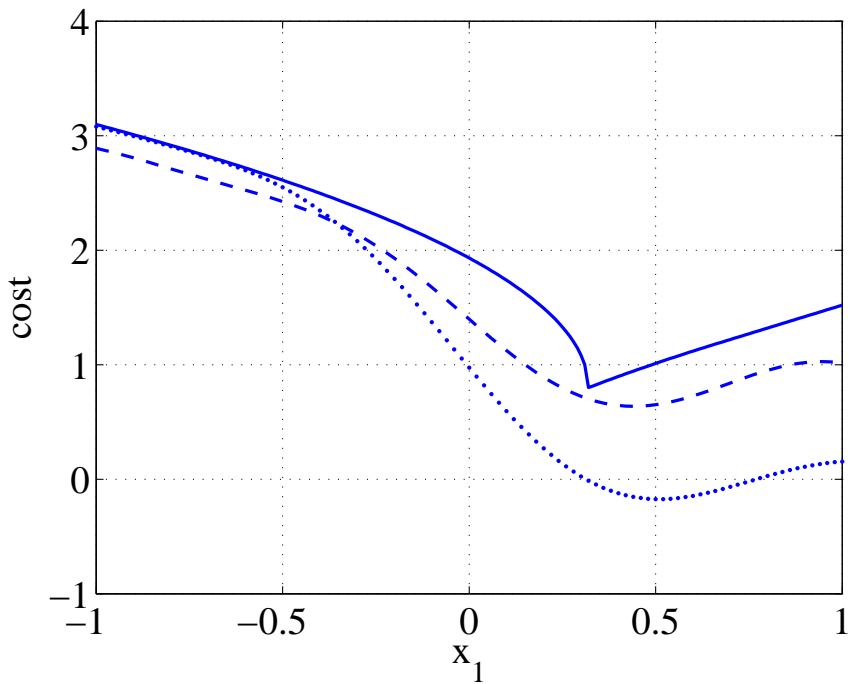} }
\caption{Example \ref{ex:dint_continued}: $T(x)$ (solid line), $\varphi(x)$ with $\SS$ being the line segment $[(-1,-0.8)-(-0.5,-0.8)]$ (dotted line) and the line segment $[(-1,-0.8)-(1,-0.8)]$ (dashed line). The graph is obtained for $x_2=-0.8$ and $r=18$.}
\label{fig:dint_continued}
\end{figure}
\end{example}

From figure \ref{fig:dint_continued} one can make the following observations:
\begin{itemize}
\item There is a trade-off between the accuracy of the approximation of the value function and the size of the set considered for the approximation;
\item As the approximating function $\varphi$ is a polynomial, it is difficult to obtain a good approximation of the optimal value function $T(x)$ at a point 
$x$ where it is not differentiable.
\end{itemize}
While the first fact is quite normal, the second deserves more attention.
The quality of the approximation and the improvement when increasing the degree of the test functions $r$ depends on the specific optimal control considered: if the value function to be approximated is smooth the moment approach performs better.

\section{Control synthesis}\label{se:control}
As already observed, when an optimal solution of the SOS problem is close to the optimum value, the (non-negative) integrand on the right hand side of equation (\ref{eq:intdiff}) takes small values along an optimal trajectory ($0$ when the HJB condition is satisfied). Therefore,  given an optimal solution $\varphi$ of the SOS relaxation (\ref{eq:sosocp}), a natural control law candidate $u(x(t))$ is a global minimizer of 
\begin{equation}\label{eq:control}
\min_{u\in U(t,x)} \!\left[\frac{\partial\varphi(x,t)}{\partial t} \!+\! \nabla_x \varphi(x,t)f(t,x,u) \!+\! h(t,x,u)\right].
\end{equation}


Suppose $g_{T_j}$ does not depend on $t$ (i.e. $g_{T_j}(x,u)$)
and define a box $\SS_{x}$ around $x$ and contained in
the set $\{x:~\exists u:~g_{T_j}(x,u)\leq 0,~ \forall j\}$.
The control algorithm we propose is the following:
\vskip 0.2cm
\noindent\framebox[\linewidth][l]{
\begin{minipage}{0.92\linewidth}\vskip 0.1cm
\begin{enumerate}
\item Set $\bar x = x(t)$.
\item Calculate the moments corresponding to a uniform 
probability measure on
$\SS_{\bar x}$.
\item Solve the moment relaxation to the OCP.
\item Apply the control obtained by minimizing (\ref{eq:control}) until
$x(t)\notin\SS_{\bar x}$.
\item Go to step 1.\end{enumerate}\vskip 0.1cm
\end{minipage}
}\vskip 0.2cm

First we show how the control strategy can be applied to the double integrator
with state and input constraints.
\begin{example}\label{ex:dint_continuedbis}
The set $U(x)$ of feasible controls is the interval $[-1,1]$ when $x_2>-1$ and $[0,1]$ when $x_2=-1$.
Indeed, when the trajectory constraint is active ($g_T(x)=x_2+1=0$), 
an admissible trajectory must be such that $\dot g_T(x)\geq 0$ and therefore $u\geq 0$. As
$f$ is affine in $u$ and the cost does not depend on $u$, the control is easily 
obtained by checking the sign of $\nabla_x \varphi(x) \begin{bmatrix}0 & 1\end{bmatrix}^T$.

Figure \ref{fig:dint_traj} displays the trajectory obtained from the initial condition $x(0)=(1,1)$
with the choice $\SS_{\bar x} = \{x:|x-{\bar x}|\leq 0.05,x_2\geq-1\}$.
The simulation stopped after $700$ iterations with a time of $3.5$ seconds to reach 
a circle of radius $0.01$; this is exactly the minimum time required to reach the origin when calculated 
using $T(\cdot)$.
\begin{figure}[htb]
\centering{
\includegraphics[width=0.8\textwidth]{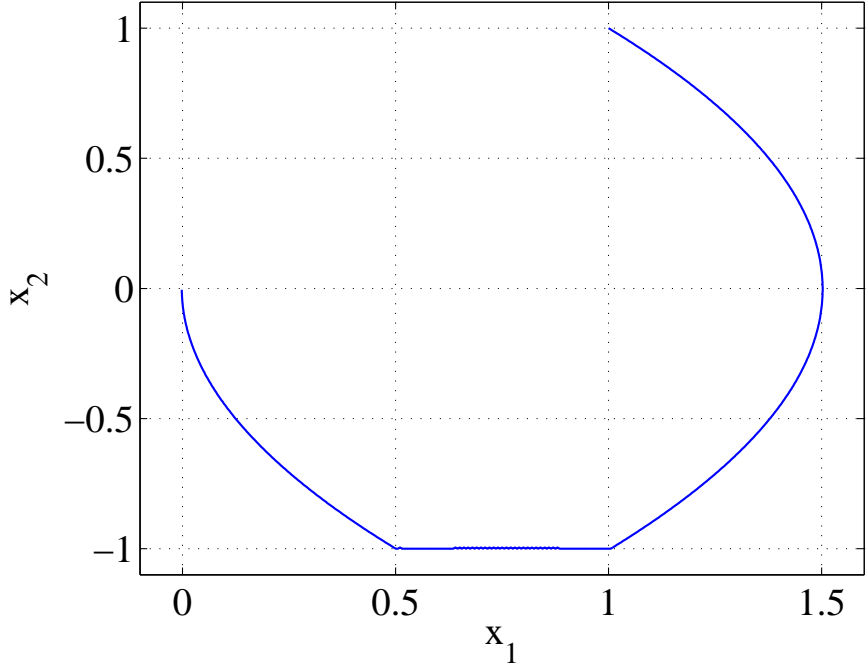} }
\caption{Example \ref{ex:dint_continuedbis}. Trajectory of the controlled system.}
\label{fig:dint_traj}
\end{figure}
\end{example}

In some cases, like in the next example, the approximate
value function can be computed once and for all at time $t=0$. It can also
be proved that the resulting control law is indeed stabilizing.
\begin{example}\label{ex:nonlinear}
Consider the nonlinear system
\begin{equation}
\begin{bmatrix} \dot x_1(t) \\ \dot x_2(t) \end{bmatrix} =
\begin{bmatrix} x_2(t) - x_1(t)^3 + x_1(t)^2 \\ u(t) \end{bmatrix}
\end{equation}
already considered in \cite{PraParRan2004}.
The objective of the optimal control problem with free terminal time 
consists of driving to the origin the initial states from
the set $(x_1,x_2)\in \SS =[-1,1]\times[-1,1]$ by minimizing the cost functional $\int_0^Th(x(t),u(t))dt$ with
\begin{equation}
h(x(t), u(t))=x_1(t)^2+x_2(t)^2+\frac{u(t)^2}{100}.
\end{equation}
From Remark \ref{re:no_t}, the approximated value function $\varphi$ 
(computed with $\mu_I$ a known initial probability distribution on $\SS$) does not depend on $t$.
As $f$ is affine in $u$ and $h$ is quadratic in $u$, the control law $u(x)$ can be obtained by the first order optimality conditions:
\begin{equation}\nonumber
u(x):=-\frac{1}{2}\nabla_x \varphi(x) \begin{bmatrix}0 & 1\end{bmatrix}^T.
\end{equation}
As this OCP has no analytic solution, we evaluate the control performance by simulating the closed-loop system considering several initial conditions and then evaluating
\begin{equation}
\mbox{gap}=\dfrac{2(UB-LB)}{UB+LB}
\end{equation}
where $LB$ is the lower bound on the cost given by the moment relaxation and $UB$ is the upper bound given integrating the cost during the simulations. The initial conditions considered are on the boundary of $\SS$ and are represented in Figure \ref{fig:nonlinear}. For such values of $x(0)$ the trajectories converge to the origin considering test functions of degree $r\geq 6$. In the next table, the maximal value of the gap for the initial conditions considered is reported for different values of $r$.
\vskip 0.2cm
\begin{center}
\begin{tabular}{|l|l|l|l|l|}
\hline degree & 6 & 8 & 10 & 12 \\
\hline gap  & 0.2275 & 0.0629 & 0.0577 & 0.0567 \\ \hline
\end{tabular}
\end{center}
\vskip 0.2cm
Observe that for $r=8$ the performance of the control is rather satisfactory. Although the gap is decreasing as expected (if optimality could be reached the gap would be $0$), the variation is really small. This is probably due to the fact that the current SDP solvers have some difficulties handling even medium size problems.
Figure \ref{fig:nonlinear} shows the trajectories obtained for $r=10$.
\begin{figure}[htb]
\centering{
\includegraphics[width=0.8\textwidth]{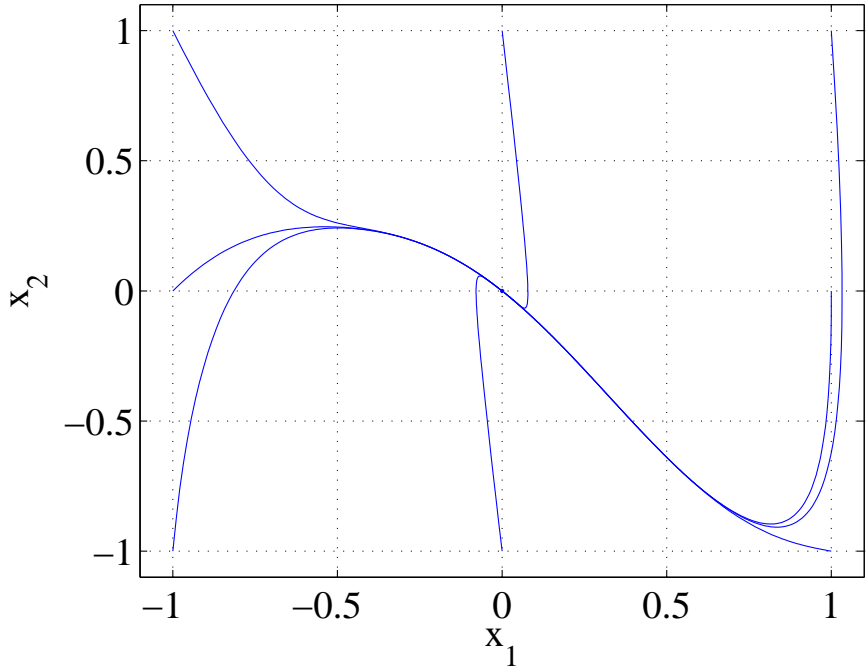} }
\caption{Example \ref{ex:nonlinear}. Trajectories obtained for $r=10$. }
\label{fig:nonlinear}
\end{figure}

Using the moment approach we can also show that the proposed control drives the state to the origin for every initial condition in $\SS$. Indeed, we consider the controlled system and the following problem: maximize $\int_0^T x_1(t)^2+x_2(t)^2 dt$ under the constraint $\CC_I=\SS$ (all the occupation measures and $T$ are undetermined). By solving this problem for $r\geq 4$, we found an upper bound on the cost. Since by linearization we can verify the origin is a local attractor, this upper bound implies that every trajectory starting from $\SS$ reaches the origin.
\end{example}

\section{Conclusions}

This paper is a follow-up to \cite{LasPriHen2005,LasHenPriTre2008}
where a sequence
of lower bounds were derived for the optimal value of a
polynomial optimal control problem (OCP), following an
occupation measure approach. In the current paper, we propose
some techniques to constructively derive a control law from
the solution of the convex linear matrix inequality (LMI)
relaxations of the OCP. So our contribution can be seen
as an extension to synthesis of the performance analysis results
of \cite{LasPriHen2005,LasHenPriTre2008}.

Generally speaking, we believe that the moment formulation
of OCP is an appealing alternative to indirect methods
based on Lyapunov or Hamilton-Jacobi-Bellman techniques.
The moment formulation deals directly with systems
trajectories. The resulting primal LMI moment problem
admit a dual LMI sum-of-squares (SOS) formulation which is,
however, instrumental to the explicit computation of a
control law. In this
context, the nice interplay between functional
analysis (measure theory) and algebraic geometry (representation
of polynomial positive on semialgebraic) may provide
constructive answers to potentially difficult
control synthesis problems.

Current limitations of the approach are as follows.

First, as we are seeking a polynomial value function (a
smooth subsolution of the Hamilton-Jacobi-Bellman
equation) that approximates the (possibly non-smooth)
optimal value function $\bar{\varphi}(t,x)$, it may
happen that precision deterioriates at points
where $\bar{\varphi}(t,x)$ is non-smooth.
Partitioning of the state-space, and/or iterative
computation of the polynomial value function in a
neighborhood moving along optimal trajectories
(like in example \ref{ex:dint_continuedbis})
could help address this issue, at the price of
an increased computational burden.

Second, we are relying on the performance of current
available general-purpose SDP solvers.
Semidefinite programming is a relatively
young research field, and the degree of maturity of SDP
solvers is far from that of, say, linear or convex
quadratic programming solvers. More specifically,
as far as we know, there is currently no
numerically stable SDP solver, and no tractable
estimate of the conditioning of an LMI problem.
For example, it is expected that the choice of a
basis to represent polynomials and moments has a
significant impact on the problem conditioning,
and hence on the numerical behavior of the solvers.

Third, the number of variables and constraints in the LMI
problems grows quickly as a function of the number of state
and input variables and the degree of the polynomial
approximation of the value function. Current general-purpose
SDP solvers can deal with a few thousands variables and
constraints, well below the dimensions of moment LMI problems
corresponding to OCPs with, say, 6 states and 2 inputs.
For these reasons, dedicated
primal-dual interior-point methods tailored to the
specific quasi-Hankel or quasi-Toeplitz structure
of moment LMI problems would be welcome.

Finally, we are currently working on a user-friendly
OCP module for GloptiPoly 3 \cite{HenLasLof2007}, that
helps formulating explicitly an OCP
as a generalized problem of moments. The user
only provides the polynomial data of the OCP, and
the module automatically generates an
approximate optimal control law. Once it is ready
and fully documented, the software will be freely
available for download from the GloptiPoly 3 webpage
\begin{center}
\tt www.laas.fr/$\sim$henrion/software/gloptipoly3
\end{center}

\section*{Acknowledgments}

This work was partly supported by project MOGA of
the French National Research Agency (ANR),
and by project No. 102/06/0652 of the Grant Agency
of the Czech Republic (GA\v CR).

\end{document}